\documentclass{article}
\pdfoutput=1
\usepackage[english]{babel}
\usepackage[utf8]{inputenc}
\usepackage{vmargin,graphicx,amsmath,amssymb,amsthm,latexsym,xcolor}
\usepackage{hyperref,fancyhdr,dsfont}
\usepackage{upgreek}
\usepackage{tikz}
\usetikzlibrary{decorations.pathreplacing}
\usetikzlibrary{decorations.pathmorphing}
\usetikzlibrary{arrows,shapes}
\usetikzlibrary{calc}
\usepackage{graphics}
\usepackage[nocompress]{cite}
\usepackage[protrusion=true,expansion,babel=true]{microtype}
\usepackage{pdfprivacy}

\usepackage{mathrsfs}
\usepackage{cancel}
\usepackage{enumerate}

\newcommand{\BE}{\begin{equation}}
\newcommand{\EE}{\end{equation}}
\newcommand{\BEn}{\begin{equation*}}
\newcommand{\EEn}{\end{equation*}}
\newcommand{\BEN}{\begin{equation*}}
\newcommand{\EEN}{\end{equation*}}
\newcommand{\BA}{\begin{array}}
\newcommand{\EA}{\end{array}}
\newcommand{\BEQNA}{\begin{eqnarray}}
\newcommand{\EEQNA}{\end{eqnarray}}
\newcommand{\BEQNAn}{\begin{eqnarray*}}
\newcommand{\EEQNAn}{\end{eqnarray*}}
\newcommand{\BEA}{\begin{eqnarray}}
\newcommand{\EEA}{\end{eqnarray}}
\newcommand{\BEAn}{\begin{eqnarray*}}
\newcommand{\EEAn}{\end{eqnarray*}}
\newcommand{\BEAN}{\begin{eqnarray*}}
\newcommand{\EEAN}{\end{eqnarray*}}
\newcommand{\BM}{\begin{multline}}
\newcommand{\EM}{\end{multline}}
\newcommand{\BMn}{\begin{multline*}}
\newcommand{\EMn}{\end{multline*}}
\newcommand{\BMN}{\begin{multline*}}
\newcommand{\EMN}{\end{multline*}}

\newcommand{\dis}{\displaystyle}

\renewcommand{\{}{\left\lbrace}
\renewcommand{\}}{\right\rbrace}
\renewcommand{\(}{\left(}
\renewcommand{\)}{\right)}

\newcommand{\R}{\mathbb{R}}
\newcommand{\N}{\mathbb{N}}

\newcommand{\h}\GHH

\newcommand{\sH}{{\rm H}}

\newcommand{\sL}{{\rm L}}

\newcommand{\HH}{\mathrm{H}}
\newcommand{\GHH}{\boldsymbol{\mathrm{H}}}
\newcommand{\WW}{\mathrm{W}}

\newcommand{\LL}{\mathrm{L}}
\newcommand{\GLL}{\boldsymbol{\mathrm{L}}}

\newcommand{\DD}{\mathrm{D}}
\newcommand{\NN}{\mathrm{N}}

\newcommand{\MA}{\mathcal{A}}

\newcommand{\MC}{\mathcal{C}}

\newcommand{\ML}{\mathcal{L}}

\newcommand{\MO}{\mathcal{O}}

\newcommand{\MP}{\mathcal{P}}

\newcommand{\Vn}{\, h \,}

\newcommand{\fonction}[5]{\begin{array}[t]{lrcl}
#1 :&#2 &\longrightarrow &#3\\
&#4& \longmapsto &#5
\end{array}}

\newcommand{\norm}[1]{\left\Vert #1 \right\Vert}
\newcommand{\abs}[1]{\left\vert #1 \right\vert}

 \usepackage{authblk}
\reversemarginpar

\usepackage{mathrsfs}
\usepackage{subfig}

\title{On the penalization by the perimeter in shape optimization applied to Dirichlet inverse obstacle problem}

\author{Fabien Caubet\footnote{Universit\'e de Pau et des Pays de l'Adour, E2S UPPA, CNRS, LMAP, UMR 5142, 64000 Pau, France. \texttt{fabien.caubet@univ-pau.fr}}, Marc Dambrine\footnote{Universit\'e de Pau et des Pays de l'Adour, E2S UPPA, CNRS, LMAP, UMR 5142, 64000 Pau, France. \texttt{marc.dambrine@univ-pau.fr}}, J\'er\'emi Dard\'e\footnote{Institut de Math\'ematiques de Toulouse, UMR 5219, Universit\'e de Toulouse, CNRS, UPS, F-31062 Toulouse Cedex 9, France. \texttt{Jeremi.Darde@math.univ-toulouse.fr}}}

\date{\textit{\today}}

\begin{document}

\maketitle
\cfoot{}
\rhead{\small \textit \leftmark}
\lhead{}
\cfoot{\thepage}
\numberwithin{equation}{section}


\begin{abstract}
This paper is devoted to the understanding of regularisation process in the shape optimization approach to the so-called Dirichlet inverse obstacle problem for elliptic operators. More precisely, we study two different regularisations of the very classical shape optimization approach consisting in minimizing a mismatched functional.
The first one is an implicit regularisation when working in the class of inclusion having a uniform $\varepsilon$-cone property, a natural class in shape optimization. As this regularity is not trivial to guarantee numerically, we discuss the regularisation by perimeter penalization. We show that this second regularisation provides a stability gain in the minimization process. 
\end{abstract}

\medskip


\medskip



\newtheorem{theor}{Theorem}[section]
\newtheorem{theorem}[theor]{Theorem}
\newtheorem{lemme}[theor]{Lemma}
\newtheorem{cor}[theor]{Corollary}
\newtheorem{prop}[theor]{Proposition}
\newtheorem{defi}[theor]{Definition}
\newtheorem{rem}[theor]{Remark}
\newtheorem{corollary}[theor]{Corollary}
\newtheorem{definition}[theor]{Definition}
\newtheorem{example}[theor]{Example}
\newtheorem{lemma}[theor]{Lemma}
\newtheorem{notation}[theor]{Notation}
\newtheorem{proposition}[theor]{Proposition}
\newtheorem{remark}[theor]{Remark}
\newtheorem{assumption}[theor]{Assumption}


\section{Introduction}


In this work, we are primarily interested in studying the effect of the penalization by the perimeter, in the context of shape optimization for 
the Dirichlet inverse obstacle problem. 
This regularisation of Tikhonov type can be seen as an analogous one to the classical inverse Cauchy problem (also called data reconstruction problem) where the adding term is linked with the norm of the unknown and permits to obtain existence and convergence results (see~\cite{CauDarGod19} and reference therein). 
Our aim is to study in details the effect of 
adding a penalization term using the perimeter when we consider the inverse problem of recovering an unknown inclusion in a domain from the knowledge of boundary measurements, in particular to quantify the stability gain, if such gain exists.
In other words the question that guided us are as follows: does the Tikhonov regularisation guarantee, and/or is-it necessary to guarantee, the existence of minimizers, depending on the class of admissible domains with which we are working? Does the addition of this term improve the stability of the numerical reconstruction?

More precisely, we consider  a bounded smooth connected open set $\Omega$ of $\R^d$ (where $d \in \N^*$ represents the dimension). In the following, $\nu$ denotes the outward normal vector to $\Omega$. Let $(g_{\NN},g_{\DD}) \in \HH^{-1/2}(\partial\Omega)\times\HH^{1/2}(\partial\Omega)$ be a non-trivial given Cauchy pair. The inverse obstacle problem we are interested in reads as follows: find $\omega$, a \textit{sufficiently smooth} domain compactly contained in~$\Omega$, and~$u\in \HH^1(\Omega\backslash \overline \omega)$, such that
\begin{equation} \label{eq_pbl_ObsDiri}
\left\lbrace \begin{array}{cccl}
- \Delta u + u &= &0& \text{ in } \Omega \backslash \overline{\omega}, \\
u  &=& g_{\DD}& \text{ on } \partial \Omega, \\
\partial_\nu u &=& g_{\NN}  & \text{ on } \partial \Omega, \\
u &=& 0&  \text{ on } \partial \omega.
\end{array}
\right.
\end{equation}
Notice that the regularity of the inclusion $\omega$ will be specified and discussed below.

The Dirichlet inverse obstacle problem is the simplest model of geometrical inverse problems, and contains the two main difficulties inherent in this class of problems: it is non-linear as we search the geometry, and it is ill-posed, as it might fail to have a solution for some pairs of Cauchy data. It is actually exponentially ill-posed, even when the geometry is known (\textit{i.e.} the problem of reconstructing $u$ from $(g_{\NN},g_{\DD})$ knowing $\omega$, see~\cite{Hadamard,BenBelgacem07,BenBelgacem05}), and therefore highly sensible to small perturbations on the data. Nevertheless, identifiability holds in the class $\omega$ Lipschitz and~$u$ continuous in a neighborhood of $\partial \omega$ (and therefore, as a consequence of elliptic regularity, for inclusions with $C^1$ boundaries): 
there exits at most one couple $(\omega,u)$ satisfying \eqref{eq_pbl_ObsDiri} (see, among others, \cite[Theorem 5.1]{Colton-Kress} or \cite[Theorem 1.1]{BouDar10}). Thus, it is reasonable to design algorithm to try to reconstruct $\omega$ and $u$ from the knowledge of a  data $(g_{\NN}, g_{\DD})$ possibly affected by some noise. From now on, we will say that the data
$(g_{\NN},g_{\DD})$ is \emph{compatible} if~\eqref{eq_pbl_ObsDiri} admits a solution, and is \emph{not compatible} otherwise.

The ill-posedness of our problem of interest is noticeable in the very poor stability estimate one can obtain for compatible data. Indeed, suppose that
we have two compatible sets of data with, for simplicity, the same Neumann condition: $(g_{\NN}, g_1)$ with corresponding solution $(\omega_1,u_1)$, and
$(g_{\NN},g_2)$ with corresponding solution $(\omega_2,u_2)$. Set $A_1 = \Omega \backslash \overline{\omega_1}$ and $A_2 = \Omega \backslash \overline{\omega_2}$.
There exists a modulus of continuity $\gamma$ with $\gamma(t)\leq C \log(t)^\eta$ where $C$ and $\eta\in (0,1]$  depend on $\Omega$ such that (see~\cite[Theorem~2-1, p.763]{Aless2000})
\begin{equation}
\label{stability:statement} 
\|g_1-g_2\|_{\sL^2(\partial\Omega)} \leq \varepsilon \implies d_{\mathcal{H}}(A_1,A_2) \leq \gamma\left( \frac{\varepsilon}{\|g_{\NN}\|_{\sH^{-1/2}(\partial\Omega)}}\right),
\end{equation}
where $d_{\mathcal{H}}$ is the so-called Hausdorff distance. Such stability result is optimal (see~\cite{Aless2000}), hence the very poor stability (logarithmic stability) of the Dirichlet inverse problem.

Numerous methods have been developed to tackle this problem (see, e.g., \cite{BouDar10,Dar12,potthast,rundell,haddar-kress}, and the references therein). In our study, we focus on the shape optimization approach (used, e.g., in~\cite{AfrDamEppKat07,CauDarGod19}), with a very standard misfit function. More precisely, let us define a set of admissible domains
$$
\MA := \left\lbrace A = \Omega \backslash \overline{\omega}, \ \theta \text{ open, } \omega \subset K, \ \omega \text{ lipschitz } \right\rbrace.
$$
Here, $K$ is a known compact set satisfying the assumption $\omega \subset K \subset \Omega$. It is classically used to prevent the open sets to be arbitrary close to the boundary of $\Omega$, which would add technical difficulties.
We consider the problem of minimizing the least-squares functional $\ML \, : \, \MA \longrightarrow \R^+$ given, for all $A \in \MA$, by
$$
\ML (A) = \frac{1}{2} \int_{\partial\Omega} \abs{u-g_{\DD}}^2 ,
$$
where $u \in \HH^1(A)$ is the solution of the following problem
\begin{equation} \label{EDPprincipale}
\left\lbrace
\begin{array}{rclcl}
- \Delta u + u & = & 0 & & \mbox{ \rm in } \, A = \Omega \backslash \overline{\omega}  ,\\
\partial_{\nu} u & = & g_{\NN} & & \mbox{ \rm on } \, \partial\Omega  ,\\
 u & = & 0 & & \mbox{ \rm on } \, \partial\omega  .
\end{array}
\right.
\end{equation}
The classical underlying idea of the approach is that if $A$ satisfies $\ML(A)=0$, then $\omega = \Omega \backslash \overline{A}$ is indeed the searched obstacle as $(\omega,u)$ satisfies \eqref{eq_pbl_ObsDiri}. Equivalently, if the Dirichlet inverse obstacle problem admits a solution, then $\ML$ admits a minimizer $A$ and $\ML(A) = 0$. Now, if the data is corrupted by some noise (which is obviously the case in practical situations), the situation gets complicated: if the corrupted data is still compatible, $\ML$ still admits 0 as minimum, but the corresponding minimizer might be distant from the correct obstacle due to the very low stability of the problem. And if the corrupted data is not compatible, the infimum of $\ML$ may be still zero but it may be that no minimizer exists. 

As expected, the problem of minimizing $\ML$ is very weakly stable, as one can obtain a logarithmic stability similar to \eqref{stability:statement} for compatible data. 
In practice, when ones considers a numerical minimization scheme for $\ML$, one does not see such poor stability. One could think the stability is too weak with respect to the numerical errors made in the scheme: when one goes in the iterations of the descent method, the boundary of the reconstructed inclusion oscillates more and more. Afraites and al have studied this in~\cite{AfrDamEppKat07} by a second order shape sensitivity analysis and showed that the Riesz representative of the shape hessian is compact: the difficulty lays at the continuous level. This clearly underlines why a regularisation is needed.

Here, we study two types of regularisation. The first one consists in enforcing more regularity on the admissibles sets $\omega$, by forcing them to satisfy an \emph{uniform $\varepsilon$-cone condition}. In other word, we replace the set $\MA$ by a more demanding set $\MA_\varepsilon$. This alone restores stability in the minimization problem, but at the same time seems not very practical from a numerical point of view. %
The second one is the standard penalization by the perimeter. It consists in replacing the function $\ML$ by its penalized version 
$$
\ML_{\eta}: A \in \MA \mapsto  \ML(A) + \eta \mathcal{P}(A),
 $$
where $\eta \geq 0$ is the regularisation parameter and  
 $\mathcal{P}(A)$ is the perimeter of $A$ (relative to~$\R^d$) given by (see, e.g.,~\cite[Definition~2.3.1]{HP})
 $$
 \MP(A)= \sup \{ \int_{A} \mathrm{div} \phi  ; \, \phi \in \boldsymbol{\mathcal{D}} (\R^d), \, \norm{\phi}_{\infty} \leq 1 \}.
 $$
 We recall that if $A$ is regular (that is $\mathcal{C}^1$), the previous definition coincides with the usual definition of perimeter: $\dis \MP(A)= \int_{\partial A} 1$ (see, e.g.,~\cite[Proposition~2.3.3]{HP}). Notice again that even if it is a very common choice of regularisation, it is not obvious that it restores well-posedness of the minimization problem in the class $\MA$. 
 It nevertheless has a positive impact in the minimization process, as in a certain sense, it improves the coerciveness of the functional at a minimizer. We quantify this improvement, which constitutes the main result of this study.
 
\paragraph{Outline of the article}
The papers is organised as follows.
We first introduce the possible optimisation problems related to $\mathcal{L}$ and $\mathcal{L}_\eta$, discussing in particular different admissible spaces of open sets. This naturally leads to several problems for which we prove existence of minimizers, and discuss the behavior of these minimizers as the different parameters of regularization go to zero.
Finally we quantify the gain of stability obtained by penalizing the functional with the perimeter, 
using the second order shape derivatives of the misfit functional.

\section{Several class of admissible sets and several regularisations}

\subsection{Four minimisation problems}

As mentioned above, a natural strategy in order to solve the inverse obstacle problem~\eqref{eq_pbl_ObsDiri} is to minimize the functional 
$$
\mathcal{L}(A) = \frac{1}{2} \int_{\partial \Omega} \vert u - g_{\DD} \vert^2 , 
$$
where $u$ is the $A$-dependant solution to~\eqref{EDPprincipale}, over the previously defined set 
$$
\mathcal{A} = \left\lbrace A = \Omega \backslash \overline{\omega},\ \omega \text{ open, } \omega  \subset K,\ \omega \text{ lipschitz } \right\rbrace,
$$
as we know that, by assumption, our target $\omega$ belongs to $\mathcal{A}$. In other words, we are interested in the minimization problem
\begin{equation} \label{eq_def_P00} \tag{$\mathcal{P}$} 
    \inf_{A \in \mathcal{A}} \ML(A).
\end{equation}
However, as said in the introduction, problem 
\eqref{eq_def_P00}
is ill-posed and highly unstable. It is therefore risky to solve this problem numerically without additional regularisation. 

A first method to regularize problem 
\eqref{eq_def_P00}
consists in restricting the set of admissible domains. 
Very classically in shape optimization, for a positive real parameter $\varepsilon > 0$, we are interested in open sets satisfying the so-called $\varepsilon$-cone property whose definition is recalled below. 

\begin{definition}
A bounded open set $\Theta$ of $\R^d$ is said to satisfy the $\varepsilon$-cone property if $$
\forall x \in \partial\Theta, \, \exists \xi_{x} \, \mbox{ unit vector such that } \,  \forall y \in \overline{\Theta} \cap \mathcal{B}(x, \varepsilon), C(y, \xi_{x}, \varepsilon) \subset \Theta,
$$
where $C(y, \xi, \varepsilon) = \{ z \in \R , (z-y,\xi) \geq \cos(\varepsilon)|z-y| \, \mbox{ and } \, 0<|z-y|<\varepsilon \}$ is the cone of vertex~$y$, of direction $\xi$ and dimension $\varepsilon$.
\end{definition}

\noindent Having this in mind, we introduce the set $\mathcal{A}_\varepsilon$ of domains having an $\varepsilon$-cone property:
$$
\mathcal{A}_\varepsilon = \left\lbrace A := \Omega \backslash \overline{\omega},\ \omega \text{ open, } \omega  \subset K,\ \omega^c \text{ has the } \varepsilon \text{-cone property} \right\rbrace.
$$
Then our first regularized problem of interest consists in minimizing the functional $\ML$ over the set~$\mathcal{A}_\varepsilon$, for a fixed positive parameter~$\varepsilon > 0$. We denote 
\eqref{eq_def_Pe0}
this second minimization problem:
\begin{equation} \label{eq_def_Pe0} \tag{$\mathcal{P}_{\varepsilon0}$}
    \inf_{A \in \mathcal{A_\varepsilon}} \ML(A).
\end{equation}
Let us make two comments on the set $\mathcal{A}_\varepsilon$. First of all, as any open set having the $\varepsilon$-cone property for some $\varepsilon>0$ is Lispchitz (see~\cite{HP}), we have $\mathcal{A}_\varepsilon \subset \mathcal{A}$ for all $\varepsilon>0$. Additionally, for all~$A \in \mathcal{A}_\varepsilon$, problem~\eqref{EDPprincipale} is well-posed. Finally, for any open set $\omega_0$ with Lipschitz boundary, there exists
$\varepsilon_0$ such that $\omega_0$ satisfies the $\varepsilon$-cone property for all $\varepsilon$ in $(0,\varepsilon_0)$ (see~\cite{HP}). Therefore, as $\Omega$ is by assumption smooth, hence Lipschitz, there exists a positive parameter $\varepsilon_*$ such that for all $\varepsilon \in (0,\varepsilon_*)$, $A \in \mathcal{A}_\varepsilon$ implies $A$ has the $\varepsilon$-cone property. From now on, the parameter $\varepsilon$ is supposed to be in the interval $(0,\varepsilon_*)$.

The second method of regularisation that we consider consists in minimising the previously defined penalized functional $\ML_\eta$ for some positive parameter $\eta >0$. Then we have two options: minimizing $\ML_\eta$ over the initial set $\mathcal{A}$, or the restricted set $\mathcal{A}_\varepsilon$, which gives two new minimization problem:
\begin{equation} \label{eq_def_P0n}
\tag{$\mathcal{P}_{0\eta}$}
    \inf_{A \in \mathcal{A}} \ML_\eta(A),
\end{equation}
and 
\begin{equation} \label{eq_def_Pen} \tag{$\mathcal{P}_{\varepsilon\eta}$}
    \inf_{A \in \mathcal{A}_\varepsilon} \ML_\eta(A).
\end{equation}

In the following, we mainly focus on problems 
\eqref{eq_def_P00}, \eqref{eq_def_Pe0} and \eqref{eq_def_Pen}. Problem \eqref{eq_def_P0n}
corresponds to a regularisation of problem 
\eqref{eq_def_P00}
by penalization of the perimeter only. It is not clear that such penalization is sufficient to ensure the well-posedness of problem 
\eqref{eq_def_P0n}. 
Indeed, in most known situations arising in shape optimization in which this perimeter penalization ensure the existence of an optimal design, the obtained minimizer is of finite perimeter but has no regularity whatsoever (see \cite{Bucur-Buttazzo}). Therefore, it seems doubtful that the penalization by the perimeter guarantees the existence of a minimzer in the class of Lispchitz domains~$\mathcal{A}$.

\subsection{Well-posedness of problems \eqref{eq_def_Pe0} and \eqref{eq_def_Pen}}

Contrary to the reformulation of the obstacle problem as the minimization problem 
\eqref{eq_def_P00}, 
there always exists a solution to the regularized minimization problems 
\eqref{eq_def_Pe0} and \eqref{eq_def_Pen}.
Let us emphasize the crucial assumption of the uniform cone property that insures the existence of a minimizer. Without it, the homogenization phenomenon may occur leading to nonconvergent minimizing sequences. 

However, a minimizer can be associated to a positive value of the reconstruction functional~$\ML$. Hence, it may not be a solution of the original inverse problem. 
This occurs, for example, if the Cauchy data are not compatible.
Additionally,  we cannot guarantee  uniqueness of the minimizer.

\begin{theorem}\label{Thm:ExistenceSansReg}
The minimisation problem \eqref{eq_def_Pe0} admits a solution denoted by $A^* \in \MA_\varepsilon$ and the penalised problem \eqref{eq_def_Pen} admits also a solution denoted by $A_{\eta}^*\in\MA$.
\end{theorem}

\begin{proof}
We focus on problem 
\eqref{eq_def_Pe0}. 
The idea is to prove that the function $\ML$ is lower semi-continuous (and even continuous in our case) for one of the Hausdorff, characteristic function, or compact sets convergence (see \cite[Chapter 2.2, Different topologies on domains]{HP}). Let $(A_{n}) \in \MA_\varepsilon$ be a sequence of domains converging to  $ A^* \in \MA_\varepsilon$ for one of these convergences. Let us prove that $ \liminf \ML(A_{n}) \geq \ML(A^*)$. Let us consider, for all $n \in \N$, the solution $u_{n} \in \HH^1(A_{n})$ of the following problem
$$
\left\lbrace
\begin{array}{rclcl}
- \Delta u_{n} + u_{n} & = & 0 & & \mbox{ \rm in } \, A_{n} = \Omega \backslash \overline{\omega_{n}}  ,\\
\partial_{\nu} u_{n} & = & g_{\NN} & & \mbox{ \rm on } \, \partial\Omega  ,\\
 u_{n} & = & 0 & & \mbox{ \rm on } \, \partial\omega_{n},
\end{array}
\right.
$$
and $u \in \HH^1(A)$ solution of the same problem with $A_n$ replaced by $A$, and $\omega_n$ replaced by $\omega$. We extend $u_n$ by zero in $\omega_n$ (resp. $u$ by zero in $\omega$), and we still denote $u_n$ (resp. $u$) this extension. By definition, both $u_n$ for all $n$ and $u$ are $\HH^1(\Omega)$ functions, and we now show that the sequence $(u_n)$ converges to $u$ in $\HH^1(\Omega)$.

Let $n\in \N$. We classically have
$$
\norm{u_{n}}_{\HH^1(\Omega)} \leq \norm{g_{\NN}}_{\HH^{-1/2}(\partial\Omega)}.
$$
Therefore, there exists $u^* \in \HH^1(\Omega)$ such that the sequence $(u_{n})$ weakly converges to $u^*$ in $\HH^1(\Omega)$. We now characterize the limit $u^*$ through a boundary value problem.
Let $\phi \in \MC_{c}^{\infty}(\Omega)$, and for all $n$ in $\N$,  $w_{n}=\phi u_{n}$. We have
$$
\left\lbrace
\begin{array}{rclcl}
\Delta w_{n} & = & \Delta \phi \, u_{n} + 2 \nabla \phi \cdot \nabla u_{n} + \phi u_{n} & & \mbox{ \rm in } \, A_{n} = \Omega \backslash \overline{\omega_{n}}  ,\\
w_{n} & = & 0 & & \mbox{ \rm on } \, \partial\Omega  ,\\
w_{n} & = & 0 & & \mbox{ \rm on } \, \partial\omega_{n}  .
\end{array}
\right.
$$
Denote $f_{n} = \Delta \phi \, u_{n} + 2 \nabla \phi \cdot \nabla u_{n} + \phi u_{n}$ and $f = \Delta \phi \, u^* + 2 \nabla \phi \cdot \nabla u^* + \phi u^*$. Since~$(u_{n})$ weakly converges to $u^*$ in $\HH^1(\Omega)$,  $(f_{n})$ converges strongly in $\HH^{-1}(\Omega)$ to $f$. From \cite[Corollary 3.2.6, Theorem 3.2.13 and Corollary 3.2.15]{HP}, we deduce that $(w_{n})$ converges strongly in $\HH^1_{0}(\Omega)$ to $w$, where $w \in \HH^1_{0}(\Omega)$ solves
$$
\left\lbrace
\begin{array}{rclcl}
\Delta w & = & f & & \mbox{ \rm in } \, A^* = \Omega \backslash \overline{\omega^*}  ,\\
w & = & 0 & & \mbox{ \rm on } \, \partial\Omega  ,\\
w & = & 0 & & \mbox{ \rm on } \, \partial\omega  .
\end{array}
\right.
$$
We notice that $w = \phi u^*$ (by uniqueness of the limit). Then, since 
$$
\Delta \phi \, u^* + 2 \nabla \phi \cdot \nabla u^* + \phi u^* = f = \Delta w = \Delta \phi \, u^* + 2 \nabla \phi \cdot \nabla u^* + \phi \Delta u^*,
$$ 
for all $\phi \in \MC_{c}^{\infty}(\Omega)$, we have $- \Delta u^* + u^* = 0$ in $ \Omega \backslash \overline{\omega^*}$. Moreover, $w = \phi u^*=u^* = 0$ on $\partial\omega^*$. 
It remains to determine the boundary condition on the fixed boundary $\partial \Omega$. To that end, we choose $\phi$ such that $\phi=1$ in $K$ and localize in the fixed domain $\Omega\setminus K$. Finally we have 
$$
\Delta ((1-\phi)u_{n}) = \Delta (1-\phi) \, u_{n} + 2 \nabla (1-\phi) \cdot \nabla u_{n} + (1-\phi) u_{n} \quad \mbox{ in } \Omega ,
$$
$(1-\phi)u_{n} = 0$ on $\partial K$, and the weak convergences $ (1-\phi)u_{n} \rightharpoonup  (1-\phi)u^* $ in $\HH^1(\Omega \backslash K)$ and $\Delta ((1-\phi)u_{n}) \rightharpoonup  \Delta((1-\phi)u^*) $ in $\LL^2(\Omega \backslash K) $. Thus, on $\partial\Omega$, we have 
$$
\partial_{\nu}((1-\phi)u_{n}) = \partial_{\nu}(u_{n}) = g_{\rm N} \rightharpoonup \partial_{\nu}((1-\phi)u^*) = \partial_{\nu}(u^*)
$$
weakly in $\HH^{-1/2}(\partial\Omega)$, and then $\partial_{\nu}u^* = g_{\rm N}$. To sum up, $u^*$ satisfies
$$
\left\lbrace
\begin{array}{rclcl}
- \Delta u^* + u^* & = & 0 & & \mbox{ \rm in } \, A^* = \Omega \backslash \overline{\omega^*}  ,\\
\partial_{\nu} u^* & = & g_{\NN} & & \mbox{ \rm on } \, \partial\Omega  ,\\
 u^* & = & 0 & & \mbox{ \rm on } \, \partial\omega^* ,
\end{array}
\right.
$$
which immediately implies $u^* = u$ by uniqueness of the solution to that problem.

Let us now prove that the weak convergence $u_{n} \rightharpoonup u$ is actually strong in $\HH^1(\Omega)$. We have $u_{n} = (1-\phi)u_{n} + \phi u_{n}$. We know that $(w_{n})=(\phi \, u_{n})$ converges strongly in $\HH^1_{0}(\Omega)$ to $w=\phi \, u$. Moreover, setting $h_{n} = (1-\phi)u_{n}$, for all $n \in \N$, and $h = (1-\phi)u$, we have
$$
\Delta h_{n} = u_n - \Delta(\phi u_n) \longrightarrow  u -  \Delta (\phi u) = \Delta h ,
$$
strongly in $\HH^{-1}(\Omega)$, and $\partial_{\nu}h_{n}=g_{\rm N} = \partial_{\nu} h$ on $\partial\Omega$. Thus, there exists a constant $c$ such that $h_n$ strongly converges to $h+c$ in $\HH^1(\Omega)$. As $h_n$ weakly converges to $h$ in $\HH^1(\Omega)$, the constant is zero.
Finally, we conclude that $u_{n}= (1-\phi)u_{n}+ \phi u_{n}$ converges strongly to $(1-\phi)u+ \phi u=u$ in $\HH^1(\Omega)$. 
Using this convergence, we obtain that $\ML$ is continuous and then $ \lim \ML(A_{n_{k}}) = \ML(A^*)$.

To conclude the proof, let $(A_{n})$ be a minimizing sequence. Then there exists a subsequence~$(A_{n_{k}})$ which converges to $A^* \in \MA$ in the sense of Hausdorff, of the characteristic functions and of the compact sets (see, e.g., \cite[Theorem~2.4.10]{HP}). Since $\ML$ is continuous, we have 
$$
\inf_{A\in\MA} \ML(A) = \lim \ML(A_{n_{k}}) = \ML(A^*) ,
$$
which proves that $A^*$ is a minimizer.

The statement for the penalized criterion follows  since the perimeter $\MP$ is lower semicontinuous for the convergence of characteristic functions (see, e.g., \cite[Proposition 2.3.7]{HP}) and since $\ML$ is continuous.

\end{proof}

\begin{remark}
Let us underline that if we consider $- \Delta u = 0$ instead of $- \Delta u + u = 0$, we obtain $\abs{u_{n}}_{\HH^1(\Omega)} \leq \norm{g_{\NN}}_{\HH^{-1/2}(\partial\Omega)}$ instead of $\norm{u_{n}}_{\HH^1(\Omega)} \leq \norm{g_{\NN}}_{\HH^{-1/2}(\partial\Omega)}$ in the previous proof. Then, in this case, we have to check that there exists a constant $c\geq 0$, independent of $n$, such that $\norm{\nabla u_{n}}_{\GLL^2(\Omega)} \geq c \norm{u_{n}}_{\HH^1(\Omega)}$ in order to conclude that there exists $u^* \in \HH^1(\Omega)$ such that $u_{n} \rightharpoonup u^*$.
\end{remark}

\subsection{Convergence results}

In this section, we prove the principal results of convergence regarding our problems of interest. To put it simply, we prove below that family of solutions of problems
\eqref{eq_def_Pen}
converges to solutions to corresponding problems 
\eqref{eq_def_Pe0}
when $\eta$ goes to zero. And if 
\eqref{eq_def_P00}
has a minimizer, the same situation occurs between 
\eqref{eq_def_Pe0}
and 
\eqref{eq_def_P00}
as $\varepsilon$ goes to zero, which leads, roughly speaking, to the diagram
$$
\eqref{eq_def_Pen} \xrightarrow{\eta \rightarrow 0} \eqref{eq_def_Pe0} \xrightarrow{\varepsilon \rightarrow 0}\eqref{eq_def_P00} .
$$
We also show that, for $\varepsilon$ small enough, we can go directly from 
\eqref{eq_def_Pen}
to 
\eqref{eq_def_P00}, that is
$$
\eqref{eq_def_Pen}\xrightarrow{\eta \rightarrow 0}\eqref{eq_def_P00}
$$
Obviously, these diagrams are simplifications, as for example most of convergences are only true up to  subsequences. The precise statements are given in the following propositions.

\begin{proposition} \label{prop_link_Pe0_Pen}
Let $\varepsilon>0$ be fixed.
For $\eta>0$, let $A_\eta$ be a  minimizer of 
\eqref{eq_def_Pen}. 
There exists a sequence $(\eta_n)$ of positive parameters converging to zero
and $A^* \in \mathcal{A}_\varepsilon$ such that $A^*$ is a minimizer of 
\eqref{eq_def_Pe0}
and $A_{\eta_n}$ converges to $A^*$ in the Hausdorff, characteristic functions and compact set topologies.
\end{proposition}

\begin{proof}
As by definition, $(A_\eta)_{\eta >0}$ is a family of elements of $\mathcal{A}_\varepsilon$, there exits a sequence of positive real numbers $(\eta_n)$ that converges to $0$ and
$A^*$ in $ \mathcal{A}_\varepsilon$ such that $A_{\eta_n}$ converges to $A^*$ in the Hausdorff, characteristic functions and compact set topologies \cite[Theorem 2.4.10]{HP}. We have already seen that $\ML(A_{\eta_n}) \xrightarrow{n \rightarrow 0} \ML(A^*)$.

Let $A$ in $\mathcal{A_\varepsilon}$. We have
$$
\ML(A^*) \xleftarrow{n \rightarrow \infty} \ML(A_{\eta_n}) \leq \ML_{\eta_n}(A_{\eta_n}) \leq \ML_{\eta_n}(A) = \ML(A) + \eta_n \mathcal{P}(A)
\xrightarrow{n\rightarrow\infty} \ML(A).
$$
The result follows.
\end{proof}

\begin{proposition} \label{prop_link_P_Pe0}
Suppose that problem 
\eqref{eq_def_P00}
admits a minimizer $A^*$. 
There exists $\varepsilon_0$ such that~$A^*$ is a minimizer of 
\eqref{eq_def_Pe0}
for all $\varepsilon$ in $(0,\varepsilon_0)$.
\end{proposition}

\begin{proof}
The minimizer $A^*$ is lipschitz by definition, hence there exists $\varepsilon_0$ such that $A^*$ belongs to $\mathcal{A}_\varepsilon$ for all $\varepsilon\in (0,\varepsilon_0)$. Let $\varepsilon$ be in $(0,\varepsilon_0)$ and $A_\varepsilon$ be a minimizer of 
\eqref{eq_def_Pe0}.
By definition, $\ML(A_\varepsilon) \leq \ML(A^*)$. But $A_\varepsilon$ belongs to $\mathcal{A}$, hence $\ML(A^*) \leq \ML(A_\varepsilon)$. The results follows.
\end{proof}

\begin{corollary}
Suppose that the inverse obstacle problem \eqref{eq_pbl_ObsDiri} admits a (necessarily unique) solution. Then 
\begin{itemize}
    \item $\mathcal{P}$ admits a unique minimizer $A^*$, and $\ML(A^*) = 0$,
    \item for $\varepsilon>0$ and $\eta>0$, let $A_{\varepsilon\eta}$ be a solution of 
    \eqref{eq_def_Pen}.
    There exists $\varepsilon_0>0$ such that
    for all~$\varepsilon \in (0,\varepsilon_0)$, there exists a sequence of positive parameters $(\eta_n)$ such that $A_{\varepsilon\eta_n}$ converges to~$A^*$
    in the Hausdorff, characteristic functions and compact set topologies.
\end{itemize}
\end{corollary}

\begin{proof}
The first point has already been proven in the introduction. Let us therefore focus on the second point: by Proposition \ref{prop_link_P_Pe0}, there exists $\varepsilon_0$ such that $A^*$ is a minimizer of
\eqref{eq_def_Pe0}
for all $\varepsilon$ in~$(0,\varepsilon_0)$.  Let $A_\varepsilon$ be another minimizer of 
\eqref{eq_def_Pe0}.
We have $\ML(A_\varepsilon) = \ML(A^*) = 0$, which immediately implies $A_\varepsilon = A^*$ by uniqueness of the solution of \eqref{eq_pbl_ObsDiri}. Then the result is a direct consequence of Proposition \ref{prop_link_Pe0_Pen}.
\end{proof}

\section{A increased stability result for the penalized problem}

As seen in the previous section, the penalization of the misfit functional by the perimeter only is a priori not sufficient to ensure the well-posedness of the minimization problem. 
At the theoretical level, we need to work in the regularized set $\mathcal A_\varepsilon$ to ensure the existence of minimizers, with or without the parameter penalization.

Nevertheless, this does not mean that the perimeter regularization is useless in our context. First of all, the set $\mathcal A_\varepsilon$ is not really practical at the numerical level, whereas the penalization by the perimeter is easily implementable. But more importantly, it turns out the regularization by the perimeter improves the stability of our problem. We intend to highlight that in the present section,  
by giving a quantitative estimate of the stability increase at a minimum.

\subsection{Preliminaries and reminders}

In order to introduce the main ideas of this section devoted to achieving an increased stability result,  let us first focus on the very classical problem of the minimization of the so-called \emph{energy functional}. The energy functional $\mathcal{E}$ of 
$\Omega$ is defined by
 $$
 \mathcal{E}(\Omega) = \inf \{ \frac{1}{2} \int_{\Omega} |\nabla u|^2  - \int_{\Omega} u ; \ u\in \sH^1_0(\Omega) \}.$$  
 The energy is homogeneous of degree $d+2$ and a volume constraint is needed to obtain existence of a minimizer of this energy. It is well known that the minimizer is then the ball and Brasco {\it et al.} have even given in a quantitative version of this inequality while studying the Faber-Krahn inequality. Precisely, they have shown in~\cite{Brasco} that there exists a dimensional constant $\sigma(d)$ such that for any open set $\Omega$ with finite measure 
 \begin{equation}
 \label{stabilite:E}
 \mathcal{E}(\Omega) \geq \mathcal{E}(B_1) + \sigma(d) \mathcal{F}(\Omega)^2,
 \end{equation}
 where $\mathcal{F}$ is the \emph{Fraenkel asymmetry} defined as 
 $$\mathcal{F}(\Omega)= \inf \{ \cfrac{|\Omega\triangle B|}{|B|};\  B \text{ ball such that }|B| =|\Omega| \},$$
 where the symbol $\triangle$ stands just here for the symmetric difference between sets. This a typical stable situation. Notice that the gap between the value at a generic domain and the value at the minimizer is express in term of a crude $L^1$ type term. This is typical of the case where no regularity assumption is made on the deformations. Let us emphasize a point already noticed by Dambrine and Lamboley in \cite{Dambrine-Lamboley}: the local stability statement is directly connected to the regularity of the class of domains. For exemple, the previous general statement can even be more precise under additional regularity assumption. If one consider a quasispherical cap $\Omega_h$, that is a domain such that its boundary if the graph over the sphere $\partial B$ of a~$\mathcal{C}^{2,\alpha}$ function $h$ 
 \begin{equation}
 \label{def:sph:cap}
 \partial\Omega_h=\{ x\in \R^d, \ x=(1+h(y))y ,\ y \in \partial B\},
 \end{equation}
 then there exists $\delta(d,\alpha)>0$ such that 
 $$\|h\|_{2,\alpha} \leq \delta(d,\alpha) \quad \Rightarrow \quad \mathcal{E}(\Omega_h)\geq \mathcal{E}(B_1) +\cfrac{1}{32d^2} \|h\|^2_{\sH^{1/2}(\partial B)}$$
as shown in \cite{Dambrine-Lamboley}. Notice that, this type of result can be the first step to derive general stability results by using the selection principle as in \cite{Brasco} or geometric constraint as in \cite{Prunier}.

Finally let us consider the regularisation of $\mathcal{E}$ by adding a the penalization term $\eta \mathcal{P}(\Omega)$, where~$\eta>0$. Since the ball is the minimizer for both the energy function $\mathcal{E}$ and the perimeter $\mathcal{P}$, thus the ball is also trivially a minimizer for the penalized energy $\mathcal{E}_\eta = \mathcal{E} + \eta \mathcal{P}$. The stability issue is then well-known thanks to the quantitative inequalities. The quantitative isoperimetric inequality has been proven by Fusco and al in \cite{Fusco} and states that there exists a dimensional constant $C(d)>0$ such that 
$$
\mathcal{P}(\Omega)\geq \mathcal{P}(B_1) +c(d) \mathcal{F}(\Omega)$$
holds for any Borel set $\Omega$ with finite measure $|B_1|$ and with a finite perimeter $P(\Omega)$ in the distributional sense of Caccioppoli-De Giorgi. Combined with 
\eqref{stabilite:E}, we obtained the following quantitative stability estimate:  
$$E_\eta(\Omega)\geq E_\eta(B_1) + (\sigma(d)+\eta c(d)) \mathcal{F}(\Omega),$$
that holds for any open set $\Omega$  with finite measure $|B_1|$ and with a finite perimeter $P(\Omega)$ in the distributional sense of Caccioppoli-De Giorgi.
In the previous statement no gain is observed while it is known from a long time by numericians that adding this penalization by the perimeter permits to regularize this optimization problem. 

The reason of this apparent contradiction lays in the class of domains under consideration. We aim to justify properly this fact at the continuous level by using the methods of shape sensitivity of order two. Using the results detailed in \cite{Dambrine-Lamboley}, we immediately get for spherical caps defined in~\eqref{def:sph:cap}: there exists positive numbers $\epsilon$ and $C(d)$  such that if $\|h\|_{2,\alpha}<\epsilon$, then 
\begin{equation*}
\label{stab:precisee:locale}
\mathcal{E}_\eta (\Omega_h)\geq \mathcal{E}_\eta (B_1) +\cfrac{1}{32d^2} \|h\|^2_{\sH^{1/2}(\partial B_1)} +\eta C(d) \|h\|^2_{\sH^{1}(\partial B_1)}.
\end{equation*}
In other words, the presence of the perimeter terms increases the local stability by passing from a $\sH^{1/2}$-stability to  a $\sH^{1}$-stability.

\subsection{Stability result concerning the inverse obstacle problem}

In order to consider the least-squares tracking objective $\mathcal{L}$, we need to consider more general shapes than balls. To that end, we first recall some basic facts from differential geometry and fix the notations. We follow in this paragraph the constructions first given in~\cite{DamPie00} and followed in \cite{Dambrine-Lamboley}. We follow step by step the general strategy developed in \cite{Dambrine-Lamboley} to prove quantitative stability results for smooth perturbations. 

In the following, we will assume that $\omega$ is a smooth (at least $\MC^{2,\alpha}$) bounded open set and we recall that $\nu$ denotes the outer unit normal field to $\partial\omega$. There is a non negative real $\rho(\Omega)$, such that the  application $T_{\partial\omega}$ defined by
$$
\fonction{ T_{\partial\omega}}{\partial\omega \times \R}{\R^d}{(M,V)}{M\ +\ V \, \nu(M)}
$$
is a local diffeomorphism from $\partial\omega\times \(-\rho(\omega),\rho(\omega) \)$ on an open tubular neighborhood of $\partial\omega$ we will denote $\mathcal{T}_{{\partial\omega}}$. This property expresses the fact that any point $x$ in 
$\mathcal{T}_{\partial\omega}$ has a unique orthogonal projection $p_{\partial\omega}(x)$ on $\partial\omega$ and that the relation 
$$x= p_{\partial\omega}(x) +  h(x) \, \nu(p_{\partial\omega}(x)) \; \text{ where } |h(x)|=\| x- p_{\partial\omega}(x)\|,$$
holds.   Notice that $h$ is uniquely defined and we can defined \textit{normal graph} on $\partial\omega$, that is $\omega_{h}$ such that
$$
\partial\omega_{h} = \{ x + h(x) \nu(x), \; x \in \partial\omega \} .
$$

Before stating and proving the following theorem, we briefly recall the definition of our main tool, that is the shape derivatives, and fix some notations. Given a domain $A$, we first define for any $V \in \WW^{1,\infty}(\mathbb{R}^d)$ such that $\mathrm{supp}(V) \subset K \subset \Omega$ and for any shape functional~$J$ the function $\mathcal{J}_{A}(V) = J\left[ (\mathrm{I}+V)(A)\right]$ on a neighborhood of $0$ in $\WW^{1,\infty}$. In particular, the first shape derivative is $\mathrm{D}J(A) \cdot V :=\mathcal{J}_{A}'(0)$, a continuous linear form on $\WW^{1,\infty}$ (the shape gradient), and the second order shape derivative is $\mathrm{D}^2J(A) \cdot[V,V]:=\mathcal{J}_{A}''(0)$, a continuous symmetric bilinear form on~$\WW^{1,\infty}$ (the shape hessian).

In the following, we assume that the perturbation direction $V$ belongs to $\MC^3(\partial\omega)$ and, taking into account the structure theorem of the shape derivatives (see~\cite[Theorem~5.9.2]{HP}), we also assume that $V=(V \cdot \nu)\nu$ on $\partial \omega$ so that  $h=V \cdot \nu$ with the previous notations.
In that context, our main result reads:

\begin{theorem}
\label{thm:stab:LS}
Let $A^*_\eta = \Omega \backslash \overline{\omega^*_\eta}$ be a domain of class $\MC^{2,\alpha}$. Then the penalised least-squares functional~$\ML_{\eta}$ is twice Fr\'echet differentiable on a neighborhood of $A^*_\eta $ for the $\WW^{1,\infty}$ norm. Moreover, if $A^*_\eta = \Omega \backslash \overline{\omega^*_\eta}$ is a critical shape for $\ML_{\eta}$ stable in the sense that
$$D^2 \ML_\eta(A^*_\eta) \cdot [V,V]>0 , \quad \forall V \in \WW^{1,\infty}(\partial\omega^*_{\eta}),$$
then there exists $\xi > 0$ and $c = c(\xi) > 0$ such that
$$
\ML_{\eta} (A_h) \geq \ML_{\eta} (A^*_\eta) + c \, \|h \|_{\HH^{1}(\partial\omega^*_\eta)}^2 ,
$$
for all $A_h = \(\mathrm{I}+h\nu \)(A^*_\eta)$ such that $\norm{h}_{\MC^{2,\alpha}(\partial\omega^*_\eta)} \leq \xi $.
\end{theorem}

In other words, the presence of the penalization by the perimeter provides a local $\sH^1$-stability at a minimum of the functional, therefore making the numerical research of a minimizer slightly more stable.

\begin{proof}
The fact that $\ML_{\eta}$ is twice Fr\'echet differentiable on a neighborhood of $A^*_\eta$ for the $\WW^{1,\infty}$ norm is classical and is based on the use of the implicit function theorem (we refer for example to~\cite[Chapter 5]{HP} or also~\cite{AfrDamEppKat07} for details). 

Let us now prove the main part of the theorem which is the stability result. 
We recall that we consider $V \in \mathrm{W}^{1,\infty}(\mathrm{R}^d)$ with $V \in \MC^{2,\alpha}(\partial\omega^*_\eta)$, $\mathrm{supp}(V) \subset K \subset \Omega$ and with $V=h \nu$ on $\partial \omega^*_\eta$ and that $A^*_\eta = \Omega \backslash \overline{\omega^*_\eta}$ is a stable critical shape for $\mathcal{L}_\eta$.

{\it Step 1: computation of the first order shape derivative.} We know that the shape derivative $u' \in \HH^1(A^*_\eta)$ solves (see, e.g., \cite[Chapter 5]{HP} or~\cite{AfrDamEppKat07})
\begin{equation*}
\left\lbrace
\begin{array}{rclcl}
- \Delta u' + u' & = & 0 & & \mbox{ \rm in } \, A^*_\eta  ,\\
\partial_{\nu} u' & = & 0 & & \mbox{ \rm on } \, \partial\Omega  ,\\
 u' & = & - \partial_{\nu} u \Vn & & \mbox{ \rm on } \, \partial\omega^*_\eta  .
\end{array}
\right.
\end{equation*}
We introduce the adjoint state $w \in \HH^1(A)$, solution to
\begin{equation*}
\left\lbrace
\begin{array}{rclcl}
- \Delta w + w & = & 0 & & \mbox{ \rm in } \, A^*_\eta  ,\\
\partial_{\nu} w & = & u-g_{\DD} & & \mbox{ \rm on } \, \partial\Omega  ,\\
 w & = & 0 & & \mbox{ \rm on } \, \partial\omega^*_\eta  .
\end{array}
\right.
\end{equation*}
We have
\begin{multline*}
\DD \ML (A^*_\eta) \cdot V 
	= \int_{\partial\Omega} u' (u-g_{\DD})
	= \int_{\partial\Omega} u' \partial_{\nu} w 
	= \int_{A^*_\eta} \Delta w \, u' + \int_{A^*_\eta} \nabla w \cdot \nabla u' - \int_{\partial\omega^*_\eta} u' \partial_{\nu} w \\
	= \int_{A^*_\eta} \Delta u' \, w + \int_{A^*_\eta} \nabla u' \cdot \nabla w - \int_{\partial\omega^*_\eta} u' \partial_{\nu} w
	= \int_{\partial\Omega} w \partial_{\nu} u' + \int_{\partial\omega^*_\eta} w \partial_{\nu} u' - \int_{\partial\omega^*_\eta} u' \partial_{\nu} w .
\end{multline*}
Thus, using the boundary conditions, we obtain
$$
\DD \ML (A^*_\eta) \cdot V = \int_{\partial\omega^*_\eta} \partial_{\nu} u \partial_{\nu} w \Vn ,
$$
and then, using the classical shape derivative of the perimeter
$$
\DD \ML_{\eta} (A^*_\eta) \cdot V = \int_{\partial \omega^*_\eta}\( \partial_{\nu} u \partial_{\nu} w + \eta \HH \) \Vn .
$$
Notice that if $A^*_\eta$ is a critical point, then
$$
\partial_{\nu} u \partial_{\nu} w = - \eta \HH \quad \mbox{ on } \; \partial \omega^*_\eta.
$$

{\it Step 2: computation of the second order shape derivative.} We know that the second order shape derivative $u'' \in \HH^1(A^*_\eta)$ solves (see, e.g., \cite[Chapter 5]{HP} or~\cite{AfrDamEppKat07})
\begin{equation*}
\left\lbrace
\begin{array}{rclcl}
- \Delta u'' + u'' & = & 0 & & \mbox{ \rm in } \, A^*_\eta  ,\\
\partial_{\nu} u'' & = & 0 & & \mbox{ \rm on } \, \partial\Omega  ,\\
 u'' & = & - (\partial_{\nu} u' - \nabla u \cdot \nabla_{\tau}\Vn)\Vn - \partial_{\nu } u V \cdot (-\nabla_{\tau} \Vn) & & \mbox{ \rm on } \, \partial\omega^*_\eta  .
\end{array}
\right.
\end{equation*}
Since $V=\Vn \nu$, we obtain that 
$$
 u''  =  - (\partial_{\nu} u' - \nabla u \cdot \nabla_{\tau}\Vn)\Vn
 \quad  \mbox{ \rm on } \, \partial\omega^*_\eta  .
$$
Moreover the shape derivative $w' \in \HH^1(A^*_\eta)$ solves
\begin{equation*}
\left\lbrace
\begin{array}{rclcl}
- \Delta w' + w' & = & 0 & & \mbox{ \rm in } \, A^*_\eta   ,\\
\partial_{\nu} w' & = & u' & & \mbox{ \rm on } \, \partial\Omega  ,\\
 w' & = & - \partial_{\nu} w \Vn & & \mbox{ \rm on } \, \partial\omega^*_\eta  .
\end{array}
\right.
\end{equation*}
We then have
\begin{multline*}
\DD^2 \ML (A^*_\eta) \cdot [V , V]
	= \int_{\partial\Omega} \( u'' (u-g_{\DD}) + \abs{u'}^2 \)
	= \int_{\partial\Omega} u'' \partial_{\nu} w + \int_{\partial\Omega} u' \partial_{\nu} w' \\
	= \int_{\partial\Omega} w \partial_{\nu} u'' + \int_{\partial\omega^*_\eta} w \partial_{\nu} u'' - \int_{\partial\omega^*_\eta} u'' \partial_{\nu} w    +    \int_{\partial\Omega} w' \partial_{\nu} u' + \int_{\partial\omega^*_\eta} w' \partial_{\nu} u' - \int_{\partial\omega} u' \partial_{\nu} w' .
\end{multline*}
Thus, using the boundary conditions  for $u''$ the second derivative of the state, we obtain
$$
\DD^2 \ML (A^*_\eta) \cdot [V , V] = \int_{\partial\omega^*_\eta} \partial_{\nu} w \( \partial_{\nu} u' - \nabla u \cdot \nabla_{\tau}\Vn \) \Vn   - \int_{\partial\omega^*_\eta} \partial_{\nu} w \partial_{\nu} u' \Vn + \int_{\partial\omega^*_\eta} \partial_{\nu} u \partial_{\nu} w' \Vn .
$$
Noticing that $u=0$ on the boundary $\partial\omega^*_\eta$ of the inclusion and hence $\nabla_{\tau} u =0 $ on $\partial\omega$ and then $ \nabla u \cdot \nabla_{\tau}  \Vn=0$ on $\partial\omega^*_\eta$, we arrive at 
$$
\DD^2 \ML (A^*_\eta) \cdot [V , V] = \int_{\partial\omega^*_\eta} \partial_{\nu} u \partial_{\nu} w' \Vn .
$$
Moreover we know that
$$
\DD^2 \MP(\Omega) \cdot V \cdot V = \int_{\partial\omega^*_\eta} \abs{\nabla_{\tau}\Vn}^2 + \int_{\partial\omega^*_\eta} \Vn^2 \( \HH^2 - \mathrm{Tr}\({}^t \nabla_{\tau}\nu \nabla_{\tau}\nu\)\) ,
$$
where $\mathrm{H}$ is the mean curvature of $\partial\omega^*_\eta$, and thus
$$
\DD^2 \ML_{\eta} (A^*_\eta) \cdot [V , V] = 
\int_{\partial\omega^*_\eta}  \partial_{\nu} u \partial_{\nu} w'  \Vn  
+ \eta \( \int_{\partial\omega^*_\eta} \abs{\nabla_{\tau}\Vn}^2 + \int_{\partial\omega} \Vn^2 \( \HH^2 - \mathrm{Tr}\({}^t \nabla_{\tau}\nu \nabla_{\tau}\nu\)\) \).
$$

{\it Step 3: the shape Hessian $\DD^2 \ML_{\eta} (A^*_\eta)$ is coercive in $H^1$.} 
We split the shape hessian as
$$
\DD^2 \ML_{\eta} (A^*_\eta) \cdot [V , V] =  \ell_{r} (\Vn,\Vn) + \ell_{m} (\Vn,\Vn),
$$
where we have set
$$
\ell_{r}(h,h) = \int_{\partial\omega^*\eta} 
\partial_{\nu} u \partial_{\nu} w'  
h
+ \eta \int_{\partial\omega^*\eta}h^2 \( \HH^2 - \mathrm{Tr}\({}^t \nabla_{\tau}\nu \nabla_{\tau}\nu\)\)
$$
and
$$
\ell_{m}(h,h) = \eta \int_{\partial\omega^*_\eta} \abs{\nabla_{\tau}h}^2 .
$$
We check easily that  $\ell_{r}$ is continuous on $\LL^{2}(\partial \omega^*_\eta)$ while $\ell_{m}$ is lower semi-continuous in $\HH^{1} (\partial\omega^*_\eta)$ and moreover satisfies the coercivity property 
$$\ell_{m}(h,h) > c_{1} \abs{h}^2_{\HH^{1}(\partial\omega^*)},$$ for all $h \in \MC^{\infty}(\partial\omega^*_\eta)$. As a consequence, we can apply Lemma 3-1 of \cite{Dambrine-Lamboley} and obtain that there exists a constant $C$ such that
$$\DD^2 \ML_{\eta} (A^*_\eta) \cdot [V , V] \geq C \abs{h}^2_{\HH^{1}(\partial\omega^*)}.$$
{\it Step 4: Improved Taylor property.} We want now to prove that there exist $\xi >0$ and a modulus of continuity $\mathfrak{m}$ such that for every domain $A_h = \Omega \setminus \overline{\omega} = \(\mathrm{I}+h \nu \)(A^*_\eta)$ with $\norm{h}_{\MC^3(\partial \omega^*_\eta)} \leq \xi$,
$$
\abs{\ML_{\eta}(A_h) - \ML_{\eta}(A^*_\eta) - \frac{1}{2} \DD^2 \ML_{\eta} (A^*_\eta) \cdot [h , h] } \leq \mathfrak{m} \(\norm{h}_{\MC^3(\partial \omega^*_\eta)}\)\norm{h}_{\HH^1(\partial \omega^*_\eta)}^2 .
$$
We split the functional into the two natural parts. Concerning the perimeter, we have (see \cite[Proposition 4.5]{Dambrine-Lamboley})
$$
\MP(A_h) = \int_{\partial\omega} 1 
	= \int_{\partial\omega^*_\eta} j_{\partial\omega}(h) 
	= \MP(A^*_\eta) + \MP^\prime(A).h +\frac{1}{2}  \MP^{\prime\prime} (A^*_\eta) \cdot [h , h]  + \MO \(\norm{h}_{\WW^{1,\infty}(\partial \omega^*_\eta)}\norm{h}_{\HH^1(\partial \omega^*_\eta)}^2\) ,
$$
where $j_{\partial\omega}(h)$ is the surfacic jacobian. 

Concerning the least-squares tracking, we follow the strategy of Dambrine and Lamboley \cite[section 4-2]{Dambrine-Lamboley}. As shown in \cite{Dambrine,Dambrine-Lamboley} the crucial step is to prove an \emph{improved continuity property of the second order shape derivative} : 
namely that there exist $\eta>0$ and a modulus of continuity $\omega$ such that for every domain $\Omega_{h}$ with $\|h\|_{2,\alpha} \leq \eta$, and all $t \in[0,1]:$
$$
\left| l^{\prime \prime}(t)-l^{\prime \prime}(0)\right| \leq \omega\left(\|h\|_{2,\alpha}\right)\|h\|_{\mathrm{H}^{s}}^{2}
$$
where $l: t \in[0,1] \mapsto \ML\left(A_{t}\right)$ for the path $\left(A_{t}\right)_{t \in[0,1]}$ connecting $A^*_\eta$ to $A_{h}$, and defined through its boundary
$$
\partial A_{t}=\{x+ th (x) \boldsymbol{\mathrm{n}}(x), x \in \partial A^*_\eta\}.
$$
We have to estimate the difference. We claim that:
$$
\left|\int_{\partial A_t} \partial_{\nu_t} u_t \partial_{\nu_t} w'_t \Vn_t-\int_{\partial A^*_\eta} \partial_{\nu} u \partial_{\nu} w' \Vn \right|\leq \omega\left(\|h\|_{2,\alpha}\right)\|h\|_{H^1/2}^2.
$$
The proof is technical but straightforward by following the various steps of \cite{Dambrine-Lamboley}. In particular, we can use directly the estimates on geometric quantities stated in \cite[Lemma 4-8]{Dambrine-Lamboley}. The control of the variations of the state $u$ and of the adjoint $w$ is a direct application of the sensitivity analysis of these functions with respect to the deformation fields. That sensitivity follows directly from the Implicit Function Theorem (see \cite[Proof of Theorem 5.7.4]{Dambrine-Lamboley},\cite{ Henry}). 

Using the Taylor formula with integral remainder:
$$
\ML\left(A_{h}\right)-\ML(A^*_\eta)=\ML^{\prime}(A^*_\eta) . h+\frac{1}{2} \ML^{\prime \prime}(A^*_\eta) .(h, h)+\int_{0}^{1}\left[l^{\prime \prime}(t)-l^{\prime \prime}(0)\right](1-t) d t
$$

{\it Step 6: conclusion.} We apply \cite[Theorem 1.1]{Dambrine-Lamboley} to prove that $A^*_\eta$ is an $\HH^1$-stable local minimum of $\ML_{\eta}$ in a $\MC^{2,\alpha}$-neighborhood, that is to say there exists $\xi > 0$ and $c = c(\xi) > 0$ such that for all 
$A = \(\mathrm{I}+h\nu \)(A^*_\eta)$ such that $\norm{h}_{\MC^{2,\alpha}(\partial\omega^*_\eta)} \leq \xi $
$$
\ML_{\eta}(A) \geq \ML_{\eta}(A^*_\eta) + c \| h \|_{\HH^{1}(\partial\omega^*_\eta)}^2.
$$
\end{proof}


\subsubsection*{Acknowledgements} This work has been partially supported by  by the Project TRECOS ANR-20-CE40-0009 funded by the
ANR (2021–2024), and the ANR LabEx CIMI (under grant ANR-11-LABX-0040) within the French State Programme “Investissements d’Avenir”.


\bibliographystyle{abbrv}
\bibliography{biblio}

\end{document}